\documentclass[11pt,oneside]{article}
\usepackage{graphicx}
\usepackage{amsmath}
\usepackage{latexsym}
\usepackage{amssymb}
\usepackage{amscd}
\usepackage{amsthm}
\usepackage{amsfonts}
\usepackage{eucal}
\usepackage[russian, english]{babel}
\usepackage{amssymb}
\usepackage[usenames]{color}

\renewcommand{\thesubsection}%
{\arabic{subsection}}

\theoremstyle{plain}
\newtheorem{theorem}{Theorem}[subsection]

\newtheorem{algorithm}{Algorithm}[subsection]
\newtheorem{definition}{Definition}[subsection]
\newtheorem{procedure}{Procedure}[subsection]

\theoremstyle{definition}
\newtheorem{example}{Example}[subsection]
\newtheorem{remark}{Remark}[subsection]
\newtheorem{problem}{Problem}[subsection]

\begin{document}

\tolerance 1000

\title{Quasigroups in cryptology}
\author{V.A. Shcherbacov}

\maketitle

\begin{abstract}

\noindent{We give a review of some known published applications of quasigroups in cryptology.}

\noindent {\bf Keywords:} cryptology, quasigroup, (r,s,t)-quasigroup, stream-cipher, secret-sharing system, zero knowledge protocol, authentication of a message, NLPN sequence, Hamming distance.
\end{abstract}

\tableofcontents

\subsection{Introduction}

Now the theory of quasigroups applications in cryptology goes through the period of rapid enough growth. Therefore any review of results in the given area of researches quite quickly becomes outdated. Here we give a re-written and supplemented  form of more early versions  \cite{SCERB_03, SCERB_03_1} of such kind of reviews. See also \cite{GLUKHOV_08, Fajar_Yuliawan}.

Almost all results obtained in the domain of quasigroups application in cryptology and coding theory till the end
of eighties years of the XX-th century are described in \cite{DK1, DK2, DK3}. In the present survey the  main
attention is devoted to the later articles in this direction.

It is possible to find basic facts on quasigroup theory in \cite{ VD, 1a, 2, HOP, LM, SCERB_03}.  Information on basic fact in cryptology can be found in many books, see, for example,
\cite{B1, BEUT, NM, NMAM}.

Cryptology is a science that consists of two parts: cryptography and cryptanalysis.  Cryptography is a science
on methods  of transformation (ciphering) of information with the purpose of this information protection from
an unlawful user. Cryptanalysis is a science on methods and ways of breaking down the ciphers \cite{DYa}.

In some sense cryptography is a "defense", i.e. this is a science on construction of new ciphers, but
cryptanalysis is an "attack", i.e. this is a science and some kind of "art", a set of methods on
breaking the ciphers. This situation is similar to situation with intelligence and contr-intelligence.

These two objects (cryptography and cryptanalysis) are very close and there does not exist a good cryptographer that does not know methods of cryptanalysis.

It is clear, that cryptology depends on level of development of society, of science and level of technology development.

We recall, a cipher is a way (a method, an algorithm) of information  transformation with the purpose of its defense. A key is some hidden part (usually, a little one) or parameter of a cipher.

Steganography is a set of means and methods of hiding the fact of sending (or passing) the information, for example, a communication or a letter. Now there exist methods of hiddenness of the fact of information sending by
usual post, by e-mail and so on. \index{steganography}

In this survey as Coding Theory (Code Theory) will be  meant a science on defense of information from  accidental errors caused by transformation and sending (passing) this information.

When sending the important and confidential information, as it seems to us, there exists a sense to use methods  of
Code Theory, Cryptology, and Steganography all together  \cite{K03}.

In cryptology one often uses the following Kerkhoff's (1835 - 1903) rule: an opponent (an unlawful user) knows all ciphering procedure (sometimes a part of plaintext or ciphertext) with exception of key. \index{Kerkhoff's
rule}

Many authors of books, devoted to cryptology divide this science (sometimes not paying attention to this fact) in two parts: before article of Diffie and Hellman \cite{DH} (so-called cryptology with non-public
(symmetric) key) and after this work (a cryptology with public or non-symmetric key). Practically namely Diffie and Hellman article opened new era in cryptology. Moreover, it is possible to apply these new approaches in
practice. \index{key} \index{key!non-public} \index{key!symmetric} \index{key!public} \index{key!non-symmetric}

Especially fast development of the second part of cryptology is connected with very fast development of Personal Computers and Nets of Personal Computers, other electronic technical devices in the end of XX-th century.  Many
new mathematical, cryptographical problems appeared in this direction and some of them are not solved. Solving of these problems have big importance for practice.

Almost all known construction of error detecting and error correcting co\-des, cryptographic algorithms and enciphering systems have made use of associative algebraic structures such as groups and fields, see, for example, \cite{MST, Dehornoy_04}.

There exists a possibility to use such non-associative structures as quasigroups and neo-fields in almost all branches of coding theory, and especially in cryptology.

Often the codes and ciphers based on non-associative systems show better possibilities than known codes and ciphers based on associative systems \cite{DK3, K97}.

Notice that in the last years the quantum code theory and quantum cryptology \cite{SOR,
EKERT, Z, Branciard_04} have been developed intensively. Quantum cryptology also use theoretical achievements of "usual" cryptology \cite{C_BENNET_84}.

Efficacy of applications of quasigroups in cryptology is based on the fact that quasigroups are "generalized permutations" of some kind and the number of quasigroups of order $n$ is larger than $n!\cdot
(n-1)!\cdot  ...  \cdot 2! \cdot 1!$ \cite{DK1}.

It is worth noting that several of the early professional cryptographers, in particular, A.A. Albert, A. Drisko, M.M.~Glukhov,  J.B. Rosser, E. Sch\"onhardt, C.I. Mendelson, R. Schaufler were connected with the development of
Quasigroup Theory. The main known "applicants" of quasigroups in cryptology were (and are) J. Denes and A.D. Keedwell \cite{D1, DK1, DK2, DK3, D00}.

Of course, one of the most effective cipher methods is to  use unknown, non-standard or very rare language. Probably the best enciphering method was (and is) to have a good agent.

\subsection{Quasigroups in "classical" cryptology}

There exist two main elementary methods when ciphering the information.

(i). Symbols in a plaintext (or in its piece (its bit))  are permuted by some law. The first known cipher of
such kind is cipher "Scital" (Sparta, 2500 years ago).

(ii). All symbols in a fixed alphabet are changed by a law on other letters of this alphabet. One of the first
ciphers of such kind was Cezar's cipher ($x\rightarrow x+3$ for any letter of Latin alphabet, for example $a\rightarrow d, b \rightarrow e$ and so on).

In many contemporary ciphers (DES, Russian GOST, Blowfish \cite{NM, DPP}) the methods (i) and (ii) are used with some modifications.

Trithemius cipher makes use of $26\times 26$ square array containing 26 letters of alphabet (assuming that the language is English) arranged in a Latin square. Different rows of this square array are used for enciphering various letters of the plaintext in a manner prescribed by the keyword or key-phrase \cite{B1,
K1}. Since a Latin square is the multiplication table of a quasigroup, this may be regarded as the earliest use of a non-associative algebraic structure in cryptology. There exists a possibility to develop this direction
using quasigroup approach, in particular, using orthogonal systems of binary or n-ary quasigroups.

R. Schaufler in his  Ph.D. dissertation  discussed the minimum amount of plaintext and corresponding ciphertext
which would be required to break the Vigenere cipher (a modification of Trithemius cipher) \cite{S1}. That is,
he considered the minimum member of entries of particular Latin square which would determine the square
completely.

Recently this problem has re-arisen as the problem of determining of so-called critical sets in Latin squares, see \cite{K6, D98, D99, DM, DHO, KEE_04}. See, also, articles, devoted to Latin trades, for example, \cite{BEAN}.

More recent enciphering systems which may be regarded as extension of Vigenere's idea are mechanical machines such as Jefferson's wheel and the M-209 Converter (used by U.S.Army until the early 1950's) and the
electronically produced stream ciphers of the present day \cite{K96, NM}.

During the second World War R.Shauffler while working for the German Cryptography service, developed a method of error detection based on the use of generalized identities (as they were later called by V.D. Belousov) in which
the check digits are calculated by means of an associative system of quasigroups (see also \cite{Damm}). He pointed out that the resulting message would be more difficult to decode by unauthorized receiver than in the
case when a single associative operation is used for calculation \cite{SCHAUFF}.

Therefore it is possible to assume that information on systems of quasigroups with generalized identities (see,
for example, works of Yu. Movsisyan \cite{MOVS} may be applied in cryptography of the present day.

\begin{definition} A bijective mapping $\varphi: g \rightarrowtail \varphi(g)$ of a finite group $(G,\cdot)$ onto
itself is called an orthomorphism if the mapping $\theta:g \rightarrowtail \theta(g)$ where $\theta(g) =
g^{-1}\varphi(g)$ is again a bijective mapping of $G$ onto itself. The orthomorphism is said to be in canonical
form if $\varphi (1) = 1$ where $1$ is the identity element of $(G,\cdot)$.
\end{definition}

A direct application of 
orthomorphisms
\ to cryptography is described in \cite{M2, M3}.


\subsection{Quasigroup-based stream ciphers} \label{Quasigroup based cryptosystem}

"Stream ciphers are an important class of encryption algorithms. They encrypt individual characters (usually binary digits) of a plaintext message one at a time, using an encryption transformation which varies with time.

By contrast, block ciphers  tend to simultaneously encrypt groups of characters of a plaintext message using a fixed encryption transformation. Stream ciphers are generally faster than block ciphers in hardware, and have less complex hardware circuitry.

They are also more appropriate, and in some cases mandatory (e.g., in some telecommunications applications), when buffering is limited or when characters must be individually processed as they are received. Because they
have limited or no error propagation, stream ciphers may also be advantageous in situations where transmission errors are highly probable" \cite{MENEZES}.
\index{cipher} \index{cipher!stream} \index{cipher!quasigroup-based}

Often for ciphering a block (a letter) $B_i$ of a plaintext the previous ciphered block $C_{i-1}$ is used. Notice that Horst Feistel was one of the first who proposed such method of encryption (Feistel net)   \cite{Feistel_73}.

In \cite{K96} (see also \cite{K97, K02}) C.~Koscielny has shown how qua\-si\-gro\-ups/ne\-o\-fi\-elds-ba\-sed
stream ciphers may be produced which are both more efficient and more secure than those based on groups/fields.



In  \cite{OS, MGS} it is  proposed to use quasigroups for secure  encoding.

A quasigroup $(Q,\cdot)$ and its $(23)$-parastrophe $(Q,\backslash)$  satisfy the following identities
$x\backslash (x\cdot y) = y, \,\, x\cdot (x\backslash y) =y$. The authors
propose to use this property of the quasigroups to construct a stream cipher.

\begin{algorithm} \label{ALG1} Let $A$ be a non-empty alphabet,  $k$ be a natural number, $u_i,
v_i \in A$, $i\in \{1,..., k\}$.  Define a quasigroup $(A,\cdot)$. It is clear that the quasigroup $(A,
\backslash)$ is defined in a unique way. Take a fixed element $l$ ($l\in A$), which  is called a leader.

Let $u_1 u_2... u_k$ be a $k$-tuple  of letters from $A$.
 The authors propose the following ciphering procedure   $v_1 = l \cdot u_1, v_{i}= v_{i-1} \cdot u_{i}$, $i= 2,..., k$.
  Therefore we obtain the following cipher-text $v_1v_2 \dots v_k$.

  The enciphering algorithm is constructed in the following way: $u_1= l \backslash v_1,
u_{i}= v_{i-1} \backslash v_{i}, i = 2,..., k.$
\end{algorithm}

The authors claim  that this cipher is resistant to the brute force attack (exhaustive search) and to  the statistical attack (in many languages some letters meet more frequently, than other ones).

\begin{example}

 Let alphabet $A$ consists from the letters $a, b, c$. Take the quasigroup $(A, \cdot)$:
\[
{\begin{array}{c|ccc}
\cdot & a & b & c  \\
\hline
a & b & c & a   \\
b & c & a & b   \\
c & a & b & c
\end{array}}
\]
Then  $(A, \backslash)$ has the following Cayley table
\[
{\begin{array}{c|ccc}
\backslash & a & b & c  \\
\hline
a & c & a & b   \\
b & b & c & a   \\
c & a & b & c
\end{array}}
\]

Let $l=a$ and open text is $u = b\,b\,c\,a\,a\,c\,b\,a$. Then the cipher text  is $v = c\,b\,b\,c\,a\,a\,c\,a$.
Applying the decoding function on $v$ we get $b\,b\,c\,a\,a\,c\,b\,a=u$.
\end{example}


Probably the cipher which is described 
here 
(Algorithm \ref{ALG1})
and its generalizations are now the most known and the most used quasigroup based stream-ciphers.

Authors \cite{OS} say that this cipher is resistant to the brute force attack and to the statistical one.

Cryptanalyses  of Algorithm \ref{ALG1} was made by M. Vojvoda \cite{VOIVODA_04}. He showed that this cipher is not resistant relatively to chosen ciphertext attack, chosen plaintext attack and ciphertext-only attack.

We give the following $3$-ary modification of Algorithm \ref{ALG1} \cite{Petrescu_07}. The possibility of such modification of Algorithm \ref{ALG1} was observed in \cite{SCERB_03}.

\begin{algorithm} \label{ALG2} Let $A$ be a non-empty alphabet,  $k$ be a natural number, $u_i,
v_i \in A$, $i\in \{1,..., k\}$.  Define a 3-ary quasigroup $(A,\beta)$. It is clear that  this quasigroup defines $(4! - 1)$ parastrophes  including $(14)$-, $(24)$- and $(34)$-parastrophe.

Take the fixed elements $l_1, l_2, l_3, l_4$ ($l_i\in A$), which  are  called  leaders.

Let $u_1 u_2... u_k$ be a $k$-tuple  of letters from $A$.
 The author proposes the following ciphering procedure   $v_1 = \beta(u_1, l_1, l_2), v_2 = \beta(u_2, l_3, l_4), v_{i}= \beta(u_i, v_{i-2}, v_{i-1}), i = 3,4,..., k-1$.
  Therefore we obtain the following cipher-text $v_1v_2 ...v_k$.

  The enciphering algorithm is constructed in the following way: $u_1= {}^{(14)}\beta (v_1, l_1, l_2),
  u_2= {}^{(14)}\beta (v_2, l_3, l_4), u_{i} = {}^{(14)}\beta (v_{i}, v_{i-2}, v_{i-1}), i = 3,4,..., k-1.$
\end{algorithm}

In \cite{Petrescu_07} also variants of Algorithm  \ref{ALG2} are given using $(24)$- and $(34)$-parastrophes of a ternary quasigroup.


Further development of Algorithm \ref{ALG1} is presented in 
\cite{GLIGOROSKY}.

\begin{definition}
Let $r$ be a positive integer. let $(Q, \ast)$ be a quasigroup and $a_j,  b_j \in  Q$. For each fixed $m \in Q$
define first the transformation $Q_m : Q^r \longrightarrow  Q^r$ by
$$
 Q_m(a_0, a_1, \dots , a_{r-1}) = (b_0, b_1, \dots ,  b_{r-1})
 \Longleftrightarrow
$$
$$  b_i  =
  \left\{
  \begin{array}{ll}
    m \ast a_0;& i =0 \\
    b_{i-1}\ast a_i; & 1 \leq i \leq  (r - 1).
\end{array}
\right.
$$
Then define ${\mathcal R}_1$ as composition of transformations of kind $Q_m$, for suitable choices of the
indexes $m$, as follows $${\mathcal R}_1(a_0,  a_1, \dots , a_{r-1}) = Q_{a_{0}} (Q_{a_1} \dots (Q_{a_{r-1}}
(a_0,  a_1, \dots ,  a_{r-1}))).$$
\end{definition}

\begin{definition} \cite{GLIGOROSKY}
 (Shapeless quasigroup) A quasigroup $(Q, \ast)$ of order $n$ is said to be shapeless if it is
non-commutative, non-associative, it does not have neither left nor right unit, it does not contain proper
subquasigroups, and there is no $k < 2n$ for which are satisfied the identities of the kinds:
\begin{equation} \label{shapeless}
\underbrace{x \ast (x \dots x \ast (x (x}_{k}\ast y)) = y;\,  y = ((y \ast \underbrace{ x)\ast \dots ) \ast x)
\ast x}_{k}
\end{equation}
\end{definition}
\index{quasigroup!shapeless}

\begin{remark}
Condition $k < 2n$ for identities (\ref{shapeless}) means that any left and right translation of quasigroup
$(Q, \ast)$ should have the order $k \geq (2n+1)$.
\end{remark}

In \cite{GLIGOROSKY} it is proposed to construct shapeless quasigroups using transversal approach
\cite{HALL_67}. Simple quasigroups without subquasigroups and with identity automorphism group are studied in
\cite{KUZ_DAN, Kepka_78, Izb_92, vs2}.

In the article \cite{danilo_08} it is proposed a  block cipher based on Algorithm \ref{ALG1}. Let $(Q,\ast)$ be
a quasigroup of finite order $2^d$. Using the operation $\ast$ authors define the following vector valued
 Boolean function (v.v.b.f.) $a \ast  b = c \Leftrightarrow \ast_{vv} (x_1, x_2, . . . , x_d, y_1, y_2,
. . . , y_d) = (z_1, z_2, . . . , z_d)$, where $x_1 . . . x_d, y_1 . . . y_d, z_1 . . . z_d$ are binary
representations of $a, b, c$ respectively.

Each element $z_i$ depends on the bits $x_1, x_2, . . . , x_d, y_1, y_2, . . . , y_d$ and is uniquely determined
by them. So, each $z_i$ can be seen as a $2d$-ary Boolean function $z_i = f_i(x_1, x_2, . . . , x_d, y_1, y_2, .
. . , y_d)$, where $f_i : \{0, 1\}^{2d}\rightarrow  \{0, 1\}$  strictly depends on, and is uniquely determined
by $\ast$.

Authors state that for every quasigroup $(Q, \ast)$ of order $2^d$ and for each bijection $Q \rightarrow \{0, 1
. . . , 2^d - 1\}$ there are a uniquely determined v.v.b.f. $\ast_{vv}$ and $d$ uniquely determined $2d$-ary
Boolean functions $f_1, f_2, . . . , f_d$ such that for each $a, b, c \in  Q$
\[
\begin{array}{l}
a \ast b = c \Leftrightarrow \ast_{vv} (x_1, . . . , x_d, y_1, . . . , y_d) = \\ (f_1(x_1, . . . , x_d, y_1, . .
. , y_d), ..., f_d(x_1, . . . , x_d, y_1, . . . , y_d)).
\end{array}
\]
 Each $k$-ary Boolean function $f(x_1, . . . ,
x_k)$ can be represented in a unique way by its algebraic normal form (ANF), i.e., as a sum of products
\[
ANF(f) = \alpha_0 + \sum_{i=1}^{k}\alpha_i x_i + \sum_{1 \leq i \leq j \leq k}^{k}\alpha_{i,j} x_i x_j + \sum_{1
\leq i \leq j \leq s \leq k}^{k}\alpha_{i,j,s} x_i x_j x_s + . . . ,
\]
 where the coefficients $\alpha_0, \alpha_{i}, \alpha_{i,j} , . . .$ are in the set $\{0, 1\}$ and the addition and
 multiplication are in the field $GF(2)$.

  The ANFs of the functions $f_i$ give  information about the complexity of the quasigroup (Q, .) via
the degrees of the Boolean functions $f_i$. The degrees of the polynomials $ANF(f_i)$ rise with the order of the
quasigroup. In general, for a randomly generated quasigroup of order $2^d$, $d \geq 4$, the degrees are higher
than 2.

\begin{definition} A quasigroup $(Q,\ast)$ of order $2^d$ is called
Multivariate Quadratic Quasigroup (MQQ) of type $Quad_{d-k}Lin_k$ if exactly $d - k$ of the polynomials $f_i$
are of degree 2 (i.e., are quadratic) and $k$ of them are of degree 1 (i.e., are linear), where $0 \leq k < d$
\cite{danilo_08}.
\end{definition}

Authors prove the following
\begin{theorem} Let $A1 = [f_{ij}]$ and $A2 = [g_{ij}]$ be two $d \times d$ matrices of linear Boolean expressions,
and let $b_1 = [u_i]$ and $b_2 = [v_i]$ be two $d\times 1$ vectors of linear or quadratic Boolean expressions.
Let the functions $f_{ij}$ and $u_i$ depend only on variables $x_1, . . . , x_d$, and let the functions $g_{ij}$
and $v_i$ depend only on variables $x_{d+1}, . . . , x_{2d}$. If $Det(A_1) = Det(A_2) = 1$ in $GF(2)$ and if
$$A_1 \cdot (x_{d+1}, . . . , x_{2d})^T + b_1 \equiv A_2 \cdot (x_1, . . . , x_d)^T + b_2$$ then the vector valued
operation $\ast_{vv}(x_1, . . . , x_{2d}) = A_1 \cdot (x_{d+1}, . . . , x_{2d})^T +b_1$ defines a quasigroup
$(Q,\ast)$ of order $2^d$ that is MQQ \cite{danilo_08}.
\end{theorem}

The authors researched the existence  of MQQ of order 8, 16 and 32.
\begin{problem}
Finding MQQs of orders $2^d$, $d \geq 6$ the authors consider as an open research problem.
\end{problem}

Authors show that the proposed cipher is resistant relatively to the chosen plain-text attack, attacks with differential
cryptanalysis, XL attack, Grobner basis attacks and some other kind of attacks.

Algebraic cryptanalysis of MQQ public key cryptosystem  is given in  \cite{Mohamed}:
"...  we present an efficient attack of the multivariate Quadratic Quasigroups (MQQ) cryptosystem. Our cryptanalysis
breaks MQQ cryptosystems by solving systems of multivariate quadratic polynomial equations using a modified version of the MutantXL algorithm".

In order to make Algorithm \ref{ALG1} more complicate and quite fast we propose the following

\begin{procedure} \label{ALG3} Let $A$ be a non-empty alphabet,  $k$ be a natural number, $u_i,
v_i \in A$, $i\in \{1,..., k\}$.  Define a system  of $n$ n-ary orthogonal operations $(A,f_i)$, $i = 1, 2, \dots, n$.
We propose  the following ciphering procedure   $v_i = f_i(u_1, u_2, \dots, u_n)$, $i = 1, 2,..., n$.
  Therefore we obtain the following cipher-text $v_1v_2 ...v_n$.

The enciphering algorithm is based on the fact that  orthogonal system of n n-ary operations
\[ \left\{
  \begin{array} {l}
    f_1(x_1, x_2, \dots , x_n) = a_1 \\
  f_2(x_1, x_2, \dots , x_n) = a_2 \\
  \dots \\
    f_n(x_1, x_2, \dots , x_n) = a_n
    \end{array}
\right.
\]
  has a unique solution for any tuple of elements $a_1, \dots , a_n$.
 \end{procedure}

Notice that we can take as a system of orthogonal $n$-ary operations a set of orthogonal  n-quasigroups \cite{Stojakovic_86, Syrbu_90, DUD_SYRB}.

Of course this choice  does not make  Procedure \ref{ALG3} more safe, but it gives a possibility to use Algorithm \ref{ALG2} and Procedure \ref{ALG3} together  on the base of the same quasigroup system.

Probably there exists a sense to use in Algorithm  \ref{ALG2} the irreducible 3-ary or 4-ary finite quasigroup \cite{AKIV_GOLD_00, AKIV_GOLD_01}.

\subsection{Some applications of quasigroup-based stream ciphers}

In \cite{OS} (see also \cite{MGS}) it is proposed to use Algorithm \ref{ALG1} for secure  encoding of file system. A survey of security mechanisms in mobile communication systems is in \cite{VOJVODA_02}.

SMS (Short Message Service) messages are sometimes used for the interchange of confidential data such as social security number, bank account number, password etc. A typing error in selecting a number when sending such a
message can have severe consequences if the message is readable to any receiver.

Most mobile operators encrypt all mobile communication data, including SMS messages. But sometimes, when encrypted, the data is readable for the operator.

Among others these needs give rise for the need to develop additional encryption for SMS messages, so that only accredited parties are able to be engaged in a communication. In \cite{HASSINEN} an  approach to this problem using
Algorithm \ref{ALG1} is described. In \cite{HASSINEN_1} differential cryptanalysis of the quasigroup cipher is
given. Definition of the encryption method is presented.

In \cite{MGS} the authors introduce a stream cipher with almost public key, based on quasigroups for defining suitable encryption and decryption. They consider the security of this method. It is shown that the key
(quasigroups) can be public and still has sufficient security. A software implementation is also given.

In \cite{KMUL} a public-key cryptosystem, using generalized quasigroup-based streamciphers is presented.  It is
shown that such a cryptosystem allows one to transmit securely both a cryptogram and a secret portion of the
enciphering key using the same insecure channel. The system is illustrated by means of a simple, but nontrivial,
example.

\subsection{Neo-fields and left neo-fields}

A left neo-field $(N, +, \cdot)$ of order $n$ consists of a set $N$ of $n$ symbols on which two binary operations
"$+$" and "$\cdot$" are defined such that $(N, +)$ is a loop, with identity element, say $0$. $(N\backslash \{0\},
\cdot )$ is a group and the operation "$\cdot$" distributes from the left over "$+$". (That is, $x\cdot(y + z) = x\cdot y + x\cdot
z$ for all $x, y, z \in  N$.) If the right distributive law also holds, the structure is called a neofield.

A left 
neofield (or neofield) whose multiplication group is $(G, \cdot)$ is said to be based on that group.
Clearly, every left neofield based on an abelian group is a neofield. Also, a neofield whose operation of addition satisfies the associative law is a field.

In \cite{DK3, DK5} some cryptological applications of neo-fields and left neo-fields are described.

\subsection{On one-way function}

A function $F: X \rightarrow Y$ is called one-way function, if the following conditions are fulfilled:

  \begin{itemize}
    \item there exists a polynomial algorithm of calculation of $F(x)$ for any $x\in X$;
    \item there does not exist a polynomial algorithm of inverting of the function $F$, i.e. there does not exist any polynomial time algorithm for solving the equation $F(x) =y$ relatively variable $x$.
\end{itemize}
\index{function!one-way}

It is proved that the problem of existence of one-way function is equivalent to well known problem of coincidence of classes P and NP.

One of  better candidates to be an one-way function is so-called function of discrete logarithms \cite{LM}.

A neofield $(N, +, \cdot)$ of order $n$ consists of a set $N$ of $n$ symbols on which two binary operations "$+$" and "$\cdot$" are defined such that $(N, +)$ is a loop with identity element, say $0$, $(N\backslash \{0\}, \cdot)$
is a group and the operation "$\cdot$" distributes from the left and right over "$+$" \cite{DK3}.

Let $(N, +, \cdot)$ be a finite Galois field or a cyclic ($(N\backslash \{0\}, \cdot)$ is a cyclic group) neofield. Then each non-zero element $u$ of the additive group or loop $(N, +)$ can be represented in the form
$u=a^{\nu}$, where $a$ is a generator of the multiplication group $(N\backslash \{0\}, \cdot)$. $\nu$ is called the discrete logarithm of $u$ with base $a$, or, sometimes, the exponent or index of $u$.

Given $\nu$ and $a$, it is easy to compute $u$ in a finite field, but, if the order of the finite field is a sufficiently large prime $p$ and also is appropriately chosen it is believed to be difficult to compute $\nu$
when $u$ (as a residue modulo $p$) and $a$ are given.

In \cite{DK3} discrete logarithms are studied over a cyclic neofield whose addition is a CI-loop.

In \cite{LM} the discrete logarithm problem for the group $RL_n$ of all row-Latin squa\-res of order $n$ is
defined (p.103) and, on pages 138 and 139, some illustrations of applications to cryptography are given.

\subsection{On hash function}

In \cite{DVORSKY_OCHODKOVA, DVORSKY_OCHODKOVA1} an approach for construction of hash function using quasigroups is described.

\begin
{definition} A function $H()$ that maps an arbitrary length message $M$ to a fixed length hash value $H(M)$ is a
OneWay Hash Function (OWHF), if it satisfies the following properties:

1. The description of $H()$ is publicly known and should not require any secret information for its operation.

2. Given $M$, it is easy to compute $H(M )$.

3. Given $H(M)$ in the rang of $H()$, it is hard to find a message $M$ for given $H(M)$, and given $M$ and
$H(M)$, it is hard to find a message $M_0 (\neq M)$ such that $H(M_0 ) = H(M)$.
\end{definition}
\index{function!hash}

\begin{definition}    A  OneWay Hash Function
 $H()$ is called  Collision Free Hash Function (CFHF), if it is hard
to find two distinct messages $M$ and $M_0$ that hash to the same result $(H(M) =
H(M_0))$\cite{DVORSKY_OCHODKOVA, DVORSKY_OCHODKOVA1}.
\end{definition}
\index{hash function!collision free}

We give construction of hashing function based on quasigroup \cite{DVORSKY_OCHODKOVA}.

\begin{definition} Let $H_Q () : Q \longrightarrow Q$ be projection defined as
\begin{equation}\label{HASH_EQUALITY}
   H_Q (q_1 q_2 \dots q_n ) =
((\dots (a \star q_1 )\star q_2\star \dots) \star q_n
\end{equation}
 Then $H_Q ()$ is said to be hash function over quasigroup $(Q;\star)$. The element $a$
is a fixed element from $Q$.
\end{definition}

\begin{example} Multiplication in the quasigroup $(Q,\star)$ is defined in the following manner: $a \star b = (a
- b)\pmod {4}$. This quasigroup has the following multiplication table:
$$\begin{array}{c|cccc}
\star & 0 & 1 & 2 & 3\\
\hline
0 & 0 & 3 & 2 & 1 \\
1 & 1 & 0 & 3 & 2 \\
2 & 2 & 1 & 0 & 3 \\
3 & 3 & 2 & 1 & 0 \\
\end{array}$$
Value of hash function is $H_2 (0013) = (((2 \star 0) \star 0) \star 1) \star 3 = 2$.
\end{example}

\begin{remark}
There exists a possibility to apply $n$-ary quasigroup approach to study hash functions of such kind. Since, in fact,  equality  (\ref{HASH_EQUALITY}) defines an $n$-ary operation.
\end{remark}

\begin{remark}
We notice, safe hash function must have at least 128-bit image, i.e. $H_Q (q_1 q_2 \dots q_n )$ must consist of at least 128-digit number \cite{NMAM}.
\end{remark}

In \cite{VOEVODE, VOIVODA_04}  hash functions, proposed  in \cite{DVORSKY_OCHODKOVA, DVORSKY_OCHODKOVA1}, are
discussed. The author shows that for some types of quasigroups these hash functions are not secure.

From \cite{MARKOV_02} we give the following summary: \lq\lq In this paper we consider two quasigroup transformations $QM1\colon A^{2m}\to A^{2m}$ and $QM2\colon A^m\to A^{2m}$, where $A$ is the carrier of a quasigroup. Based on these transformations we show that different kinds of hash functions can be designed with suitable security.\rq\rq

Further development of quasigroup based on hash function is reflected in \cite{Snasel}.

In \cite{SATTI} on Algorithm \ref{ALG1} based on encrypter that has good scrambling properties is proposed.

\subsection{CI-quasigroups and cryptology}

In \cite{DK3, GOLOMB_PATENT} some applications of CI-quasigroups in cryptology with non-sym\-met\-ric key are
described.

\begin{definition} Suppose that there exists a permutation $J$ of the elements of a quasigroup $(Q,\circ)$
such that,  for all $x, y \in Q$ $$J^r(x\circ y)\circ J^s x = J^t y,$$ where $r,s,t$ are integers. Then
$(Q,\circ)$ is called an $(r,s,t)$-inverse quasigroup (\cite{KS}).
\end{definition}
In the special case when $r=t=0$, $s=1$, we have a definition of CI-quasigroup.

\begin{example} \label{EXCI} A CI-quasigroup can be used to provide a one-time pad for key exchange (without the intervention of a key distributing centre) \cite{DK3,K}.

The sender S, using a physical random number generator (see \cite{K97} on random number generator based on quasigroups), selects an arbitrary element $c^{(u)}$ of the CI-quasigroup $(Q,\circ)$ and sends both  $c^{(u)}$ and enciphered key (message) $c^{(u)}\circ m$. The receiver R uses this knowledge of the algorithm for obtaining $J c^{(u)}= c^{(u+1)}$ from $c^{(u)}$ and hence he computes $(c^{(u)}\circ m) \circ c^{(u+1)} = m.$
\end{example}

\begin{example}
We can propose the following application of rst-inverse quasigroups in situation similar to situation described
in Example \ref{EXCI}. It is possible to re-write definitive equality of rst-inverse quasigroup in the following
manner $J^r(J^k u\circ m)\circ J^{s+k} u = J^t m.$

Then the schema of the previous example can be re-written in the following manner. The sender S selects an arbitrary element $J^k u$ of the rst-quasigroup $(Q,\circ)$ and sends both $J^k u$ and enciphered key
(message) $J^r(J^k u\circ m)$. The receiver $R$ uses this knowledge of the algorithm for obtaining $J^{k+s}(u)$ from $J^{k}(u)$ and hence he computes $J^r(J^k u\circ m)\circ J^{s+k} u = J^t m$ and after this he computes the
message $m$. Of course this example can be modified.
\end{example}

\begin{example} \cite{DK3}. \label{EX8}
Take a CI-quasigroup with a long inverse cycle $(c\,c^{\prime}\, c^{\prime \prime} \dots $ $ c^{t-1})$ of length $t$. Suppose that all the users $U_i$  $(i=1,2, \dots)$ are provided with apparatus (for example, a chip card)
which will compute $a\circ b$ for any given $a, b \in Q$. We assume that only the key distributing centre has a
knowledge of the long inverse cycle which serves as a look-up table for keys.

Each user $U_i$ has a public key $u_i \in Q$ and a private key $J u_i$, both supplied in advance by the key
distributing centre. User $U_s$ wishes to send a message $m$ to user $U_t$. He uses $U_t$'s public key $u_t$ to
compute $u_t \circ m$ and sends that to $U_t$. $U_t$ computes $(u_t \circ m)\circ J u_t = m$.
\end{example}

\begin{remark}
It is not very difficult to understand that opponent which knows the permutation $J$ may decipher a message
encrypted by this method.
\end{remark}

\begin{remark}
There exists a possibility  to generalize  Example \ref{EX8} using some 
$m$-inverse quasigroups \cite{ks_02}, or  $(r,s,t)$-inverse quasigroups \cite{KS, ks3}, else  $(\alpha, \beta, \gamma)$-inverse quasigroups \cite{KEED_SCERB}.
\end{remark}

\subsection{Critical sets and secret sharing systems}

\begin{definition}  A critical set $C$ in a Latin square $L$ of order $n$ is a set $C = \{(i; j;
k) \mid i, j, k \in \{1,2, \dots, n\}\}$ with the following two properties:

(1) $L$ is the only Latin square of order $n$ which has symbols $k$ in cell $(i, j)$ for each $(i; j; k)\in C$;

(2) no proper subset of $C$ has property (1) \cite{LM}.
\end{definition}
\index{critical set} \index{Latin square!critical set of}

A critical set is called minimal if it is a critical set of smallest possible cardinality for $L$. In other
words a critical set is a partial Latin square which is uniquely completable to a Latin square of order $n$.

\index{secret sharing scheme}
 If the scheme has $k$ participants, a $(t,k)$-secret sharing scheme is a system
where $k$ pieces of information called shares or shadows of a secret key $K$ are distributed so that each
participant has a share such that

(1) the key $K$ can be reconstructed from knowledge of any $t$ or more shares;

(2) the key $K$ cannot be reconstructed from knowledge of fewer than $t$ shares.

Such systems were first studied in 1979. Simmons \cite{S9} surveyed various secret sharing schemes. Secret sharing schemes based on critical sets in Latin squares are studied in  \cite{CDS}. We note, critical sets of
Latin squares give rise to the possibilities to construct secret-sharing systems.

Critical sets of Latin squares were studied in sufficiently big number of articles. We survey results from some of these articles. In \cite{DLR}  the spectrum of critical sets in Latin squares of order $2\sp n$ is studied.
The paper \cite{DH} gives constructive proofs that critical sets exist for all sizes between $[n^2/4]$ and $[(n^2- n)/2]$, with the exception of size $n^2/4+ 1$ for even values of $n$.

For Latin squares of order $n$, the size of a smallest critical set is  denoted by ${\rm scs}(n)$ in
\cite{Cavenagh_07}. The main result of \cite{Cavenagh_07}  is that ${\rm scs}(n)\ge n\lfloor {1\over2}(\log
n)^{1/3}\rfloor$ for all positive integers $n$.

In \cite{Horak_07} the authors show that any critical set in a Latin square of order $n\geq 7$ must have at
least $\lfloor{ {7n-\sqrt{n} -20}\over{2}} \rfloor $ empty cells. See, also, \cite{Horak_02}.

 The paper \cite{D99} contains lists of (a) theorems on the possible sizes of critical
sets in Latin squares of order less than 11, (b) publications, where these theorems are proved, (c) concrete
examples of such type of critical sets. In \cite{DM} an algorithm for writing any Latin interchange as a sum of
intercalates is corrected.

In \cite{Hamalainen_07} the author proposes a greedy algorithm to find critical sets in Latin squares. He applies this algorithm to Latin squares which are abelian 2-groups to find new critical sets in these Latin squares. The critical sets have the nice property that they all intersect some $2 \times 2$ Latin subsquare in a unique element so that it is easy to show the criticality.

In \cite{BEAN_05} the author gives an example of a critical set of size 121 in the elementary abelian 2-group of order 16.

In \cite{Mojdeh_07} critical sets of  symmetric Latin squares are studied. Therefore the authors require all elements in their critical sets and uniquely completable partial Latin squares to lie on or above the main
diagonal. For $n>2$, a general procedure is given for writing down a uniquely completable partial symmetric $2n\times2n$ Latin square $L'_{2n}$ containing $n^2-n+2$ entries, of which $2n - 2$ are identical and lie on the
main diagonal.

Paper \cite{D98} presents a solution to the interesting combinatorial problem of finding a minimal number of elements in a given Latin square of odd order $n$ by which one may restore the initial form of this square. In particular, it is proved that in every cyclic Latin square of odd order $n$ the minimal number of elements equals to $n(n-1)/2$.

Surveys on critical sets of Latin squares are given in  \cite{K6, KEE_04}. See, also, \cite{KEED_07}.

The concept of Latin trades is closely connected with the concept of critical set in Latin squares. Let $T$ be a partial Latin square and $L$ be a Latin square with $T \subseteq L$. We say that $T$ is a Latin trade if there
exists a partial Latin square $T'$ with $T'\cap T=\emptyset $ such that $(L \setminus T)\cup T'$ is a Latin square. Information on Latin trades is in \cite{Cavenagh_05}.

\begin{remark} See also Introduction for other application of critical sets of Latin squares in cryptology.
\end{remark}

"For a given
triple of permutations $T = (\alpha, \beta, \gamma)$ the set of all Latin squares $L$ such that $T$ is its autotopy is denoted by $LS(T)$. The cardinality of $LS(T)$ is denoted by $\Delta(T)$. Specifically, the
computation of $\Delta(T)$ for any triple $T$  is at the moment an open
problem having relevance in secret sharing schemes related to Latin squares" \cite{Falcon_06, Falcon_07}.

\subsection{Secret sharing systems and other algebraic systems}

Some secret-sharing systems are pointed in \cite{DK2}. One of such systems is the Reed-Solomon code over a
Galois field $GF[q]$ with generating matrix $C(a_{ij})$ of size $k\times (q-1)$, $k\leq q-1$. The determinant
formed by any $k$ columns of $G$ is a non-zero element of $GF[q]$. The Hamming distance $d$ of this code is
maximal $(d = q-k)$ and any $k$ from $q-1$ keys unlock the secret.

In \cite{GB_SECR} an approach to some Reed-Solomon codes as a some kind of orthogonal systems of n-ary operations is developed.

In \cite{ BEL_QRS_09}  general approach to construction of secret sharing systems using some kinds of orthogonal systems of n-ary operations is given. Transformations of  orthogonal systems of n-ary operations are studied in \cite{BEL_09}.

We give the summary from \cite{FITINA_07} : "We investigate subsets of critical sets of some Youden squares in the context of secret-sharing schemes. A subset $\mathcal C$ of a Youden square is called a critical set if $\mathcal C$ can be uniquely completed to a Youden square but no proper subset of $\mathcal C$ has a unique completion to a Youden square."

"That part of a Youden square $Y$ which is inaccessible to subsets of a critical set $\mathcal C$ of $Y$, called the strongbox of $\mathcal C$, may be thought to contain secret information. We study the size of the secret. J.
R. Seberry and A. P. Street \cite{Seberry_00} have shown how strongboxes may be used in hierarchical and compartmentalized secret-sharing schemes."

\subsection{Row-Latin squares based cryptosystems}

A possible application in cryptology of Latin power sets is proposed in \cite{DP}. \index{power set}

In \cite{D00} an encrypting device is described, based on row-Latin squares with maximal period equal to the
Mangoldt function.

In our opinion big perspectives has an application of row-Latin squares in various branches of contemporary cryptology ("neo-cryp\-to\-lo\-gy").

In 
\cite{LM} it is proposed to use: 1) row-Latin squares to generate an open key; 2) a conventional system for transmission of a message that is the form of a Latin square; 3) row-Latin square analogue of the RSA system; 4) procedure of digital signature based on row-Latin squares.

\begin{example}

Let $$L= \begin{array}{cccc}
  2 & 3 & 4 & 1 \\
  4 & 1 & 3 & 2 \\
  3 & 2 & 4 & 1 \\
  4 & 3 & 1 & 2 \\
\end{array}$$

Then
$$L^7 = \begin{array}{cccc}
  4 & 1 & 2 & 3 \\
  4 & 1 & 2 & 3 \\
  3 & 2 & 4 & 1 \\
  3 & 4 & 2 & 1 \\
\end{array}$$

$$L^3 = \begin{array}{cccc}
  4 & 1 & 2 & 3 \\
  1 & 2 & 3 & 4 \\
  1 & 2 & 3 & 4 \\
  3 & 4 & 2 & 1 \\
\end{array}$$

Then $$L^{21} = \begin{array}{cccc}
  2 & 3 & 4 & 1 \\
  1 & 2 & 3 & 4 \\
  1 & 2 & 3 & 4 \\
  4 & 3 & 1 & 2 \\
\end{array}$$

\noindent is a common key for a user $A$ with the key $L^3$ and a user $B$ with the key $L^7$.
\end{example}

A public-key cryptosystem, using generalized quasigroup-based streamciphers, as it has been noticed earlier, is presented in \cite{KMUL}.

\subsection{NLPN sequences over GF[q]}

Non-binary pseudo-random sequences over GF[q] of length $q^m-1$ called PN sequences have been known for a long time \cite{GOL}. PN  sequences over a finite field GF[q] are unsuitable directly for cryptology because of
their strong linear structure \cite{K97}. Usually PN sequences are defined over a finite field and often an irreducible polynomial for their generation is used. \index{NLPN sequence}

In article \cite{K97} definition of PN sequence  was generalized with the purpose to use these sequences in cryptology.

We notice, in some sense ciphering is making a \lq\lq pseudo-random sequence" from a plaintext, and cryptanalysis is a science how to reduce a check of all possible variants (cases) by deciphering of some ciphertext.

These new sequences were called  NLPN-sequences (non-linear pseudo-noise sequences). C. Kos\-ciel\-ny proposed
the following method for construction of NLPN-sequences.

Let $\overrightarrow{a}$ be a PN sequence of length $q^m-1$ over GF[q], $q>2$, 
i.e.
\[
\overrightarrow{a} = a_0a_1 \dots a_{q^m - 2}. \] Let $\overrightarrow{a}^i$ be its cyclic  $i$ places shifted to the right. For example \[\overrightarrow{a}^1 = a_1 \dots a_{q^m - 2}a_0.\]
 Let $Q = (SQ,\cdot)$ be a quasigroup of order $q$ defined on the set of
elements of the field GF[q].

Then $\overrightarrow{b} = \overrightarrow{a}\cdot \overrightarrow{a}^i$,
$\overrightarrow{c} = \overrightarrow{a}^i \cdot \overrightarrow{a}$, where $b_j = a_j\cdot a_j^i$, $c_j =
a_j^i\cdot a_j$ for any suitable value of index $j$ ($j\in \{1, 2, \dots, q^m-1\}$) are called NLPN sequences 
\cite{K97}.

NLPN sequences have much more randomness than PN sequences. As notice C. Koscielny the method of construction of
NLPN sequences is especially convenient for fast software encryption. It is proposed to use NLPN sequences by
generation of keys. See also \cite{KL}.

\subsection{Authentication of a message}

By authentication of message we mean that it is made possible for a receiver of a message to verify that the
message has not been modified in transit, so that it is not possible for an interceptor to substitute a false
message for a legitimate one.

By identification of a message we mean that it is made  possible for the receiver of a message to ascertain its
origin, so that it is not possible for an intruder to masquerade as someone else.

By non-repudiation we mean that a sender should not be able later to deny falsely that he had sent a message.

In \cite{DK3} some quasigroup approaches to problems of identification of a message, problem of non-repudiation
of a message, production of dynamic password and to digital fingerprinting are discussed. See also \cite{C}.

In \cite{DK5} authors suggested a new authentication scheme based on quasigroups (Latin squares). See also
\cite{DK2, DK3, DDO}

In \cite{SS} several cryptosystems based on quasigroups upon various combinatorial objects such as orthogonal
Latin squares and frequency squares, block designs, and room squares are considered.

\begin{definition}
Let $2\le t<k<v$. A generalized $S(t,k,v)$ Steiner system is a finite block design $(T,{\cal B})$ such that (1)
$\vert T\vert=v$; (2) ${\cal B}={\cal B}'\cup{\cal B}''$, where any $B'\in{\cal B}'$, called a maximal block,
has $k$ points and $2\le\vert B''\vert<k$ for any $B''\in{\cal B}''$, called a small block; (3) for any
$B''\in{\cal B}''$ there exists a $B'\in{\cal B}'$ such that $B''\subseteq B'$; (4) every subset of $T$ with $t$
elements not belonging to the same $B''\in{\cal B}''$ is contained in exactly one maximal block.
\end{definition}
\index{Steiner system}

In \cite{MZ} (see also \cite{EM}) an application of  generalized $S(t,k,v)$ Steiner systems in cryptology is
proposed, namely, it is introduced a new authentication scheme based on the generalized Steiner systems, and the
properties of such scheme are studied in the generalized affine planes.

\subsection{Zero knowledge protocol}

In \cite{RIVEST} Rivest  introduced All-Or-Nothing (AON) encryption mode in order to devise means to make
brute-force search more difficult, by appropriately pre-processing a message before encrypting it. The method is
general, but it was initially discussed for block-cipher encryption, using fixed-length blocks. \index{zero
knowledge protocol}

It is an unkeyed transformation, mapping a sequence of input blocks $(x_1,$ $ x_2,$ $ \dots ,$ $ x_s)$  to a
sequence of output blocks $(y_1, y_2, \dots, y_t)$ having the following properties:

 Having all blocks $(y_1, y_2, \dots, y_t)$ it is easy to compute $(x_1, x_2, \dots , x_s)$.

If any output block $y_j$ is missing, then it is computationally infeasible to obtain any information about any input block  $x_J$.

 The main idea is to preserve a small-length key (e.g.
64-bit) for the main encryption that can be handled by special hardware with not enough processing power or memory. This gives the method a strong advantage, since we can have strong encryption for devices that have
minimum performance.

Several transformation methods have been proposed in the literature for AON. In the article \cite{MAB} it is proposed a special transform which is based on the use of a quasigroup (it is used in algorithm \ref{ALG1}).

In \cite{DD} it is proposed to use isotopy of quasigroups in zero knowledge protocol.

Assume the users $(u_1 ,u_2 ,...,u_k )$ form a network. The user $u_i$ has public-key $L_{u_i }$, $L_{u_i }^{ '}$ (denotes two isotopic Latin squares of order $n$) and secret-key $I_{u_i}$ (denotes the isotopism of $L_{u_i}$ upon
$L_{u_i }^{'}$). The user $u_i $ wants to prove identity for $u_j $ but he doesn't want to reveal the secret-key
(zero-knowledge proof).
\\

1. $u_i $ randomly permutes $L_{u_i } $ to produce another Latin square H.

2. $u_i $ sends H to $u_j $.

3. $u_j $ asks $u_i $ either to:

\hspace*{5mm}a. prove that H and $L_{u_i }^{ '} $ are isotopic,

\hspace*{5mm}b. prove that H and $L_{u_i } $ are isotopic.

4. $u_i $ complies. He either

\hspace*{5mm}a. proves that H and $L_{u_i }^{'} $  are isotopic,

\hspace*{5mm}b. proves that H and $L_{u_i } $ are isotopic.

5. $u_i $ and $u_j $ repeat steps 1. through 4. $n$ times.
\\

\begin{remark} In the last procedure it is possible to use isotopy of n-ary groupoids.
\end{remark}

\subsection{Hamming distance between quasigroups}

The following question is very important by construction of quasigroup based cryptosystems: how big is the distance between different binary or n-ary quasigroups? Information on Hamming distance between quasigroup operation is in
the articles \cite{DR1, DR2, DR3, DR4, DR5, DR6, PV}.

We recall, if $\alpha$ and $\beta$ are two $n$-ary operations on a finite set $\Omega$, then the Hamming distance of $\alpha$ and $\beta$ is defined by
$$\text{dist}(\alpha,\beta)=|\{(u_1,\dots,u_n)\in\Omega^n:\alpha(u_1,\dots,u_n)\ne\beta(u_1,\dots,u_n)\}|.$$
\index{Hamming distance} \index{operations!Hamming distance}

The author in \cite{DR1} discusses Hamming distances of algebraic objects with binary operations. He also explains how the distance set of two quasigroups yields a 2-complex, and points out a connection with dissections of equilateral triangles.

 For a fixed group $(G, \circ)$, $\delta(G, \circ)$ is defined to be the minimum of all such distances for
$(G, \star)$ not equal to $(G, \circ)$ and $\nu(G, \circ)$ the minimum for $(G, \star)$ not isomorphic to $(G,
\circ)$.

In \cite{DR4} it is proved that $\delta(G, \circ)$ is $6n-18$ if $n$ is odd, $6n-20$ if $(G, \circ)$ is dihedral of twice odd order and $6n-24$ otherwise for any group $(G, \circ)$ of order greater than 50.  In \cite{PV} it
is shown that $\delta(G, \circ)=6p-18$ for $n=p$, a prime, and $p>7$.

In the article \cite{DR3} there are listed a number of group  orders for which the distance is less than the value suggested by the above theorems. New results obtained in this direction are in \cite{DR6}.

\subsection{Generation of quasigroups for cryptographical needs}

Important cryptographical problem is a generation of  \,  "big" quasigroups which it is possible to keep easily in a compact form in computer memory. It is clear that for this aims the most suitable  is a way to keep a little base and some procedures of obtaining a necessary element.

Therefore we should have easily generated objects (cyclic group, abelian group, group),  fast and complicate methods of their  transformation (parastrophy, isotopy, isostrophy, crossed isotopy \cite{Shaposhnikov_99},  homotopy, generalized isotopy), their glue and blowing (direct product, semi-direct product, wreath product \cite{KM}, crossed product, generalized crossed product).
For these aims various linear quasigroups (especially $n$-ary quasigrous) are quite suitable 
\cite{2, MARS, SCERB_08}.

In \cite{NOSOV_06} the boolean function is proposed to use by construction of $n$-ary and binary quasigroups.

A method of generating a practically unlimited number of quasigroups of an  arbitrary (theoretically) order using the computer algebra system Maple 7 is presented in \cite{K02}.

This problem is crucial to cryptography and its solution permits to implement practical quasigroup-based
endomorphic cryptosystems.

In this article \cite{K02} it is proposed to use isotopy of quasigroups and direct products of quasigroups. If we start
from  class of finite groups, then, using these ways, it is possible to obtain only class of quasigroups that
are isotopic to  groups. We notice, there exists many quasigroups (especially of  large order) that are not
isotopic to a group. Therefore for construction of quasigroups that are not isotopic to groups probably better to use the concept of gisotopy \cite{MS05, SCERB_08}.

\subsection{Conclusion remarks}

In many cases in cryptography it is possible to change associative systems by non-associative ones and practically in any case this change gives in some sense better results than use of associative systems.
Quasigroups in spite of their simplicity, have various applications in cryptology. Many new cryptographical
algorithms can be formed on the basis of quasigroups.

\textbf{Acknowledgment.} The author thanks  the Supreme Council for Science and Technological Development of Republic of Moldova (grant 08.820.08.08 RF)  for financial support.

\addcontentsline{toc}{subsection}{\protect{References}}

\vspace{2mm}
\begin{center}
\begin{parbox}{118mm}{\footnotesize
V.A. Shcherbacov,  \hfill
Received November 4, 2009

\vspace{3mm}

Institute of Mathematics and Computer Science

Academy of Sciences of Moldova

Academiei 5, Chi\c sin\u au MD-2028 Moldova

E--mail: $scerb@math.md$

}

\end{parbox}
\end{center}

\end{document}